\documentclass[11pt]{article}
\usepackage{geometry}                
\usepackage{graphicx}
\usepackage{amssymb}
\usepackage{epstopdf}
\usepackage{comment}
\usepackage{mathtools}
\usepackage{float}
\usepackage{pgfplots}
\usepackage{amsthm}
\usepackage{caption}
\usepackage{float}
\usepackage[colorlinks=true,linkcolor=black,citecolor=black,urlcolor=black]{hyperref}
\usepackage[english]{babel}
\usepackage[mathscr]{eucal}
\usepackage{amsmath}
\usepackage{array}
\usepackage{pdflscape}
\usepackage{enumerate}
\usepackage{rotating}
\theoremstyle{definition}
\newtheorem{rmk}{Remark}

\DeclareMathOperator{\ep}{\epsilon}

\newcommand{\G}{\mathbf{G}}
\newcommand{\A}{\mathbf{A}}
\renewcommand{\L}{\mathbf{L}}

\title{
  On the Numerical Solution of Fourth-Order Linear Two-Point Boundary Value Problems}
\author{   William Leeb and Vladimir Rokhlin}
\date{}

\begin{document}
\maketitle
\begin{abstract}
This paper introduces a fast and numerically stable algorithm for the solution of fourth-order linear boundary value problems on an interval. This type of equation arises in a variety of settings in physics and signal processing. Our method reformulates the equation as a collection of second-kind integral equations defined on local subdomains. Each such equation can be stably discretized and solved. The boundary values of these local solutions are matched by solving a banded linear system. The method of deferred corrections is then used to increase the accuracy of the scheme. Deferred corrections requires applying the integral operator to a function on the entire domain, for which we provide an algorithm with linear cost. We illustrate the performance of our method on several numerical examples.
\end{abstract}

\section{Introduction}
This paper describes the numerical solution to differential equations of the form
\begin{align} \label{eq:main_prob}
    \sum_{j=0}^{4}a_j(x) \frac{d^j \phi}{d x^j}(x) = f(x)
\end{align}
for $x$ in an interval $[a,b]$, with specified boundary conditions
\begin{align}
    \phi(a) &= \alpha_{l,0} \label{bdry1}\\
    \phi(b) &= \alpha_{r,0} \label{bdry2}\\
    \phi^\prime(a) &= \alpha_{l,1} \label{bdry3}\\
    \phi^\prime(b) &= \alpha_{r,1} \label{bdry4}
\end{align}
and given coefficients $a_j(x)$ and right hand side $f(x)$.

Fourth-order equations of this kind arise in a variety of physical problems. In the small-bending regime, the shape of a beam under an external force is described as the solution to an equation of the form 
\begin{align}
    \frac{d^2}{d x^2} \left( c(x) \frac{d^2 \phi}{d x^2}(x) \right) = f(x)
\end{align}
where $c(x)$ is the stiffness of the beam, and $f(x)$ is an external force applied to the beam \cite{gere-mechanics-1984}. One-dimensional equations also arise from separation-of-variables for higher-dimensional problems, such as vibration of plates \cite{courant-methods-1937}. A notable fourth-order operator on the unbounded domain $[0,\infty)$ is defined by:
\begin{align}
    -\frac{d^2}{dx^2}\left(x^2 \frac{d^2 \phi}{dx^2}(x)\right) 
                  + (a^2 + b^2)\frac{d}{dx} \left(x^2 \frac{d \phi}{dx}(x)\right) 
                  - (a^2 b^2 x^2 - 2a^2)\phi(x).
\end{align}
This operator arises in many applications because it commutes with the truncated Laplace transform composed with its adjoint \cite{grunbaum1, grunbaum2}. It is useful when working with the family of decaying exponential functions, and has been the subject of recent investigation \cite{roy-thesis,lederman-laplace-i,lederman-laplace-ii, lederman2017lower}. 


One standard approach to solving fourth-order boundary value problems of the form \eqref{eq:main_prob}--\eqref{bdry4} is to employ a finite difference or finite element scheme. To obtain the solution on $m$ nodes, such a  method involves solving an $O(m)$-by-$O(m)$ banded system of linear equations. While the solution can be obtained with asymptotic cost $O(m)$, the discretizations used introduce a condition number of size $O(m^4)$ \cite{strang-fem-2008}. The resulting loss of accuracy in the solution is entirely due to the choice of discretization, and is not a result of the conditioning inherent to the problem \eqref{eq:main_prob}--\eqref{bdry4}. Other methods have also been proposed for this class of problems, including the use of the spectral method \cite{olver2013fast}, the spectral-Galerkin method \cite{shen1995efficient}, and the Sinc-Galerkin method \cite{smith1991sinc}.

Most directly relevant to our present work are those algorithms based on reformulating the differential equation \eqref{eq:main_prob}--\eqref{bdry4} as an integral equation, such as the method contained in \cite{driscoll2010automatic}.
%
%
The integral equation approach centers on expressing the solution $\phi$ in the form 
\begin{align}
    \phi(x) = \int_{a}^{b} G_0(x,t) \sigma(t) dt + \psi_\alpha(x),
\end{align}
where $G_0(x,t)$ is the Green's function for the biharmonic equation $\phi^{(4)} = f$ with zero boundary values; $\psi_\alpha$ is a third degree polynomial with the desired boundary values; and $\sigma = \phi^{(4)}$ is the new function to be solved for. The key observation is that the function $\sigma$ can be expressed as the solution to a second-kind integral equation. While direct discretizations of the differential equation \eqref{eq:main_prob}--\eqref{bdry4} are ill-conditioned, second-kind integral equations can be stably discretized. More precisely, the values of $\sigma$ on a grid of points is expressible as the solution to a linear system whose condition number does not markedly exceed the condition number of the original continuous problem.

The challenge with using second-kind integral equations is that while their discretized linear systems are well-conditioned, they are dense; consequently, a naive solver will have cubic cost $O(m^3)$. It has been observed that for many physical problems, these dense linear systems can nevertheless be solved in linear or nearly linear time. This observation underpins a variety of methods devised for the solution of second-order two-point boundary value problems \cite{starr-bvp,greengard-two-point,lee-stiff}, as well as the class of fast multipole methods \cite{greengard-fast}.

In this paper, we compute $\sigma$ using a more direct algorithm than previous methods for solving second-kind integral equations. At a high level, our approach consists of three steps. First, we use the integral equation form of the problem to produce $m$ local solutions of the equation \eqref{eq:main_prob} with homogeneous boundary values. Next, we solve a banded linear system to match the boundary values of the local solutions. While this step introduces extraneous loss of accuracy similar to a finite element scheme, it is corrected by the use of deferred corrections, in which we recursively solve for the residual solutions on the entire interval $[a,b]$. Computing the right hand side of the residual equation requires applying the integral operators defined by the biharmonic Green's function and its derivatives to an arbitrary function on the entire interval $[a,b]$, for which we provide a linear time algorithm.



The asymptotic CPU time of our algorithm is $O(m  n^3  \log(1/\epsilon))$, where $\epsilon$ is machine precision, $n$ is a user-selected integer (the number of points for each local solution), and $m$ is the number of local solutions, or discretization nodes. The complexity of the algorithm is linear in the number of discretization nodes, but nevertheless achieves full accuracy. In this sense, our algorithm realizes the advantages of finite element methods' small CPU time while maintaining the numerical accuracy afforded by second-kind integral equations.

The rest of the paper is structured as follows. In Section \ref{sec:prelim}, we review the mathematical and numerical tools we will be using in our algorithm. In Section \ref{sec:algorithm}, we describe in detail the algorithm for solving \eqref{eq:main_prob}--\eqref{bdry4}. In Section \ref{sec:experiments}, we provide the results of numerical experiments.

\section{Mathematical and numerical preliminaries} \label{sec:prelim}

In this section, we review the mathematical and numerical tools that we will be using throughout the paper. In particular, we will show how to express the function $\phi$ which solves \eqref{eq:main_prob}--\eqref{bdry4} in terms of the solution to a second-kind integral equation; review the properties of Gaussian quadrature, which we will use to discretize this integral equation; and review the method of deferred corrections.

\subsection{Rescaling the problem domain} \label{sec:rescale}

It will be convenient to rescale the problem \eqref{eq:main_prob}--\eqref{bdry4} to convert the interval $[a,b]$ into $[-1,1]$. Define:
\begin{align}
    \tilde{\phi}(y) &= \phi\left(\frac{(b-a)(y+1)}{2}+a\right) \\
    \tilde{a}_j(y) &= \left( \frac{2}{b-a} \bigg)^j 
                          a_j \bigg( \frac{(b-a)(y+1)}{2} + a \right)  \\
    \tilde{f}(y) &=f \left( \frac{ (b-a) (y+1) }{2 } + a \right).
\end{align}
The equations \eqref{eq:main_prob}--\eqref{bdry4} for $\phi$ are then equivalent to the equation
\begin{align} \label{eq:rescale}
    \sum_{i=0}^{4} \tilde{a}_j(y)   \frac{d^j \tilde{\phi}}{dy^j}(y)
    = \tilde{f}(y)
\end{align}
with modified boundary values
\begin{align}
    \tilde{\phi}(- 1) &= \alpha_{l,0} \\
    \tilde{\phi}(1) &= \alpha_{r,0} \\
    \tilde{\phi}^\prime(-1) &= \left(\frac{b-a}{2}\right)\alpha_{l,1} \\
    \tilde{\phi}^\prime(1) &= \left(\frac{b-a}{2}\right)\alpha_{r,1} .
\end{align}
After solving for $\tilde{\phi}$ and its derivatives, we can perform the inverse change of variables to arrive at the solution $\phi$ and its derivatives, as follows:
\begin{align} \label{eq:inverse_rescale}
    \phi^{(j)}(x) = \left(\frac{2}{b-a}\right)^j
                      \tilde{\phi}^{(j)} \left( \frac{2(x-a)}{b-a} - 1 \right).
\end{align}

\subsection{The biharmonic Green's function on $[-1,1]$}

We will make central use of the Green's function $G_0(x,t)$ for the biharmonic equation 
$\phi^{(4)} = f$ with homogeneous boundary conditions ($\alpha_{l,0} = \alpha_{r,0} = \alpha_{l,1} = \alpha_{r,1} = 0$ in \eqref{bdry1}--\eqref{bdry4}), on the interval $[-1,1]$. The biharmonic Green's function is given by the formula:
\begin{align} \label{eq:Green}
    G_0(x,t) = 
    \begin{cases}
        (1-t)^2(1+x)^2(1+2t-2x-tx) / 24, \,\,\, \text{ if } t>x\\
        (1-x)^2(1+t)^2(1+2x-2t-tx) / 24, \,\,\,\text{ if } t<x\\
    \end{cases}
\end{align}
It can be checked by direct calculation that $G_0(x,t)$ satisfies the defining properties of the Green's function (see \cite{courant-methods-1937}). Specifically, the following properties hold:
\begin{align}
    \frac{\partial^4 G}{\partial x^4}(x,t) = 0, \quad t \in [-1,1],
\end{align}
\begin{align}
    G_0(-1,t) = G_0(1,t) = \frac{\partial G_0}{\partial x} (-1,t) 
                         = \frac{\partial G_0 }{\partial x}(1,t) = 0,
    \quad  t \in [-1,1],
\end{align}
\begin{align}
    \lim_{t \to x^+}\frac{\partial^j G}{\partial x^j} (x,t)
    =     \lim_{t \to x^-}\frac{\partial^j G}{\partial x^j} (x,t),
    \quad  x  \in [-1,1], \quad 0 \le j \le 2,
\end{align}
and
\begin{align}
    \lim_{h \to 0} \left[\frac{\partial^3 G}{\partial x^3}(t+h,t) 
                               - \frac{\partial^3 G}{\partial x^3}(t-h,t)\right] = 1, 
                          \quad t \in [-1,1].
\end{align}

We will use the notation $G_j$ to denote the $j^{th}$ partial derivative of the Green's function, $0 \le j \le 3$; that is, we define:
\begin{align}
   G_j(x,t) = \frac{\partial^j G_0}{ \partial x^j} (x,t) .
\end{align}
We will use the boldface letter $\G_j$ to denote the corresponding integral operator. In this notation, for a function $f$ on $[-1,1]$ we will write:
\begin{align}
    (\G_j f)(x) = \int_{-1}^{1} G_j(x,t)f(t) d t 
                = \int_{-1}^{1} \frac{\partial^j G_0}{\partial t^j}(x,t) f(t) d t.
\end{align}

The significance of the Green's function and its derivatives is that the solution to the equation $\phi^{(4)} = f$ with zero boundary conditions is the function $\phi = \G_0 f$, with derivatives $\phi^{(j)} = \G_j f$, $1 \le j \le 3$. In Section \ref{sec:fredholm}, we will use the biharmonic Green's function to reformulate the problem \eqref{eq:main_prob}--\eqref{bdry4} as a second-kind integral equation.

\subsection{Integral form of the boundary value problem \eqref{eq:main_prob}--\eqref{bdry4}} \label{sec:fredholm}
In this section, we will use the biharmonic Green's function to reformulate the differential equation \eqref{eq:main_prob}--\eqref{bdry4} as a second-kind integral equation. For comprehensive background on second-kind integrals equations, see \cite{riesz-functional, courant-methods-1937, byron-mathematics-2012}.  We let $\L$ denote the differential operator on the left side of \eqref{eq:main_prob}; that is,
\begin{align}
    (\L \phi)(x)
    = \sum_{j=0}^{4} a_j(x) \frac{d^j \phi}{dx^j}(x)
\end{align}
The equation \eqref{eq:main_prob} can then be written in the more compact form
\begin{math} 
    \L \phi = f.
\end{math}
%

Any four-times differentiable function with vanishing values and first derivatives on $[-1,1]$ can be written in the form $\G_0 \sigma$, for some function $\sigma$. Consequently, we can express the solution $\phi$ to \eqref{eq:main_prob}--\eqref{bdry4} as being of the form $\phi = \G_0 \sigma + \psi_{\alpha}$, where $\sigma = \phi^{(4)}$, and where $\psi_\alpha$ is a function satisfying the boundary conditions \eqref{bdry1}--\eqref{bdry4} and $\psi_\alpha^{(4)} = 0$ on $[-1,1]$.

It can be directly checked that the unique function $\psi_\alpha$ on $[-1,1]$ satisfying $\psi^{(4)} = 0$ and the boundary conditions \eqref{bdry1}--\eqref{bdry4} (with $a=-1$ and $b=1$) is given by the following formula:
\begin{align} \label{eq:psi_alpha}
    \psi_\alpha(x) = \alpha_{l,0} \psi_{l,0}(x) 
                  + \alpha_{r,0} \psi_r(x) + \alpha_{l,1} \psi_{l,1}(x) 
                  + \alpha_{r,1}\psi_{r,1}(x)
\end{align}
where the functions $\psi_{l,0}$, $\psi_{r,0}$, $\psi_{l,1}$ and $\psi_{r,1}$ are defined by the formulas:
\begin{align}
    \psi_{l,0}(x) &= (1-x)^2 (2+x) /4   \label{eq:psi_l0}\\
    \psi_{r,0}(x) &= (1+x)^2 (2-x) /4   \label{eq:psi_r0}\\
    \psi_{l,1}(x) &= (1-x)^2 (x+1) /4   \label{eq:psi_l1} \\
    \psi_{r,1}(x) &= (1+x)^2 (x-1) /4   \label{eq:psi_r1}.
\end{align}

The function $\sigma = \phi^{(4)}$ is expressible as the solution to a second-kind integral equation, as we now show. The equation $\L (\G_0 \sigma + \psi_\alpha) = f$
can be written equivalently as
\begin{align} \label{eq:int_eq}
    \L \G_0 \sigma =  f - \L \psi_\alpha.
\end{align}
Because $\frac{\partial^4 G}{\partial x^4} (x,t) = \delta(x-t)$, the operator $\L \G_0$ is of the form
\begin{align} 
\label{LG0}
(\L\G_0\sigma)(x) 
= a_{4}(x) \sigma(x)
     +  \sum_{j=0}^{3} a_j(x) \int_{-1}^{1}  G_j(x,t) \sigma(t) dt.
\end{align}

After solving \eqref{eq:int_eq} for $\sigma$, the solution $\phi$ to \eqref{eq:main_prob}--\eqref{bdry4} and its derivatives are given by $\phi^{(j)} = \G_j \sigma + \psi_\alpha^{(j)}$, for $0 \le j \le 3$, and $\phi^{(4)} = \sigma$. The reformulation of the differential equation as a second-kind integral equation is beneficial, as the latter can be stably discretized.


\begin{rmk}
The integral equation formulation of \eqref{eq:main_prob}--\eqref{bdry4} makes use of the biharmonic Green's function $G_0$. However, in principle this can be replaced by the Green's function $G_0$ for any equation $\L_0 \phi = f$ with a fourth-order differential operator $\L_0$.
%
%
In this general setting, $G_0$ is known as the \emph{background} Green's function \cite{greengard-two-point,lee-stiff}. We choose to work with the biharmonic equation because it is so analytically tractable. In our experience, the numerical performance is unlikely to depend substantially on the choice of background Green's function, as was observed for second-order equations in the numerical experiments of \cite{greengard-two-point}.
\end{rmk}

\subsection{Legendre polynomials and Gaussian quadrature} \label{sec:legendre}
The $n^{th}$ Legendre polynomial $P_n$ is defined for $x \in [-1,1]$ by 
\begin{align}
P_0(x) = 1; \quad
P_n(x) = \frac{1}{2^n n!} \frac{d^n}{dx^n} (x^2-1)^n, \quad
n \ge 1.
\end{align}
As is well known \cite{abramowitz}, the $n$ roots $y_1,\dots,y_n$ of $P_n$ lie in 
$[-1,1]$. Together with a suitable choice of weights $w_i > 0$, they can be used to 
evaluate the integral on $[-1,1]$ of any polynomial $f$ of degree less than or equal 
to $2n+1$ via the formula
\begin{align}
\int_{-1}^{1} f(x) dx = \sum_{i=1}^{n} f(y_i) w_i.
\end{align}

By rescaling the interval $[-1,1]$ to $[a,b]$, we obtain the quadrature
\begin{align}
\label{eq:gauss_quadr}
\int_{a}^{b} f(x) dx = \sum_{i=1}^{n} f\bigg( \frac{b-a}{2} (y_i + 1)  + a \bigg) 
                                                               \frac{b-a}{2} w_i.
\end{align}

The rescaled nodes $ \frac{b-a}{2} (y_i + 1)  + a $ are the roots of the polynomials 
$P_n(2(x-a)/(b-a) - 1)$ on $[a,b]$, which form an orthogonal basis for $L^2([a,b])$. 
The nodes $y_i$ and weights $w_i$ on $[-1,1]$ can be computed at cost $O(n)$ 
\cite{rokhlin1}.

For an arbitrary $(2n+2)$-times continuously differentiable function $f$ on $[a,b]$ (not 
necessarily a degree $2n+1$ polynomial), the same quadrature formula 
\eqref{eq:gauss_quadr} has error
\begin{align}
\bigg|\int_{a}^{b} f(x) dx 
    - \sum_{i=1}^{n} f\bigg( \frac{b-a}{2} (y_i + 1)  + a \bigg) 
                     & \frac{b-a}{2} w_i \bigg|   \nonumber \\
 &\le \frac{|f^{(2n+2)}(\xi)|}{(2n+2)!} \int_{a}^{b} \prod_{i=0}^{n} (x-x_i)^2 dx
\end{align}
for some $\xi \in [a,b]$; for a proof, see Theorem 7.4.5 in \cite{dahlquist}. It follows 
immediately that the error can be bounded above by
\begin{align} \label{leg_err}
\frac{(b-a)^{2n+3}}{(2n+2)!} \sup_{x \in [a,b]} |f^{(2n+2)}(x)|.
\end{align}

We also note that the Legendre polynomials $P_n$ satisfy the three-term recurrence 
\cite{dahlquist,gradshteyn}
\begin{align}
P_{n+1}(x) = \frac{2n+1}{n+1}x P_n(x) - \frac{n}{n+1}P_{n-1}(x), \quad n\ge1
\end{align}
Using this recurrence, for any value of $x \in [-1,1]$ all values $P_0(x),\dots,P_n(x)$ 
can be computed in $O(n)$ floating-point operations.

There is a one-to-one correspondence between the values of $f$ on the $n$ nodes $y_1,\dots y_n$ and the first $n$ Legendre coefficients of $f$. More precisely, given the values of $f$ on $n$ Gaussian nodes $y_i$ of $[-1,1]$, the first $n$ coefficients of the Legendre expansion of $f$ can be computed at cost $O(n^2)$, by applying the matrix $[P_i(y_j)w_j]_{1\le i,j \le n}$ to the vector $(f(y_1),\dots,f(y_n))^\top$. Conversely, from the first $n$ Legendre coefficients of $f$ we can compute the values of $f$ on the $n$ Gaussian nodes at cost $O(n^2)$, by applying the inverse matrix to the vector of coefficients.

\subsection{Deferred corrections}
\label{sec:corrections}
In this section, we review the technique of deferred corrections for improving the accuracy of the solution to a linear system $\A \phi = f$. Informally, the method proceeds by solving a sequence of equations for the residuals of the previous solution. The setting where this method may be applied is when there already exists an inaccurate method for inverting $\A$ but a highly accurate method for applying $\A$ to an arbitrary vector. The algorithm is as follows:

\begin{enumerate}

\item
Produce an initial solution $\hat{\phi}$.

\item
Apply $\A$ to $\hat{\phi}$ and compute the residual right hand side $\Delta = f - \A \hat{\phi}$.

\item
Compute a vector $\hat{r}$ that solves the residual system $\A r = \Delta$ for the residual $r = \phi - \hat{\phi}$.
%
%
Update the solution: $\hat{\phi} \leftarrow \hat{\phi} + \hat{r}$.

\item
Repeat steps 2 and 3 until convergence.

\end{enumerate}

Algorithms based on this method have been employed for discretized second-kind integral operators $\A$ in, for example, \cite{dutt2000spectral,glaser2009solvers,kushnir2012stiff}. Its convergence properties for certain differential equations have been studied \cite{hansen2011deferred}.


At each step of the algorithm, the current solution $\hat{\phi}$ has a residual vector $r = \phi - \hat{\phi}$ that satisfies the residual linear system:
\begin{align}
    \A r = \A \phi - \A \hat{\phi} = f - \A \hat{\phi}.
\end{align}

%
%
Since $\A \hat{\phi}$ can be computed accurately, the right hand side $f - \A \hat{\phi}$ can be computed accurately as well (though see Remark \ref{rmk:itref} below). Consequently, the residual can be obtained to some relative accuracy $\epsilon$, or in other words, we can produce a vector $\hat{r}$ with $\|\hat{r} - r\| \le \epsilon\|r\|$. Then the updated solution $\hat{\phi} + \hat{r}$ satisfies:
\begin{align}
   \|(\hat{\phi} + \hat{r}) - \phi\| = \|\hat{r} - (\phi - \hat{\phi})  \|
   = \| \hat{r} - r\| \le \epsilon \|r\| = \epsilon \|\phi - \hat{\phi}\|.
\end{align}
In other words, the error of $\hat{\phi} + \hat{r}$ is smaller than the error of $\hat{\phi}$ by a factor of $\epsilon$. If this factor is gained after every iteration, then only $\lceil \log(\ep^*) / \log(\ep) \rceil$ iterations are required to achieve machine precision $\ep^*$. 

\begin{rmk} \label{rmk:itref}
We explain one subtlety with the method of deferred corrections. In order to solve for the residual $r$ to relative error $\epsilon$, the right hand side $\Delta = f - \A \hat{\phi}$ must be computed to relative error at most $\epsilon$. Even if $\A \hat{\phi}$ is computed to full machine precision, if $\phi \approx \phi$ then $f \approx \A \hat{\phi}$, and the difference $f - \A \hat{\phi}$ may have  large relative error. More precisely, when  $\hat{\phi}$ is close enough to $\phi$ so that $\| f - \A \hat{\phi}\| \le \delta \| f\|$, then the difference $f - \A \hat{\phi}$ will be computed to relative error $\ep^* / \delta$, due to loss of digits. So long as $\ep^*/\delta \le \epsilon$, the next iteration of the algorithm will decrease the error by a factor of $\ep$. However, when $\hat{\phi}$ is close enough to $\phi$ that  $\delta < \ep^*/\ep$, the next iteration will only increase the accuracy by a factor of $\ep^* / \delta$.

\end{rmk}

\section{The algorithm for solving \eqref{eq:main_prob}--\eqref{bdry4}} \label{sec:algorithm}

In this section, we provide a detailed description of the algorithm for solving \eqref{eq:main_prob}--\eqref{bdry4}. The algorithm has three main components. First, the interval $[a,b]$ is broken into $m$ subintervals $[x_{i}, x_{i+1}], 1 \le i \le m$, with user-prescribed endpoints $x_i$ satisfying $x_1=a$ and $x_{m+1}=b$. These subintervals do not need to be the same length. A solution with zero boundary values is produced on each subinterval, employing the second-kind integral equation described in Section \ref{sec:fredholm}. In addition, four linearly independent solutions to the homogeneous equation are also produced on each subinterval.

Second, we form a linear combination of the solutions on each subinterval. The coefficients of this linear combination are chosen so that the resulting functions on adjacent intervals have matching boundary values. Solving for these coefficients is inexpensive but ill-conditioned, resulting in a solution that does not achieve machine precision.

Third, we remedy the ill-conditioning introduced in the second step by deferred corrections. To do so we must accurately apply the second-kind integral operator $\L \G_0$ to the right hand side of the equation; concretely, this entails integrating the right hand side against the functions $G_j$. We describe an algorithm for doing so with asymptotic cost $O(m)$. This step is then iterated until full machine accuracy is achieved.

The algorithm can be summarized as follows:
\begin{enumerate}
\item \label{step:solve}
Solve the problem with zero boundary values on each subinterval $[x_i,x_{i+1}]$.

\item \label{step:adjust}
Adjust the boundary values on each subinterval $[x_i,x_{i+1}]$ so that the solutions match at the endpoints.

\item \label{step:defer}
Compute the residual right hand side on the full domain $[a,b]$ and apply deferred corrections.

\end{enumerate}

Steps \ref{step:solve}--\ref{step:defer} will produce the solution $\phi$ sampled at $n$ Gaussian nodes of each subinterval $[x_i,x_{i+1}]$. These $n$ nodes can be converted into the first $n$ Legendre coefficients of the solution $\phi$ on each subinterval $[x_{i},x_{i+1}]$. The solution can then be evaluated at any point in $[a,b]$ by evaluating this expansion on the subinterval $[x_i, x_{i+1}]$ containing that point, as described in Section \ref{sec:legendre}.

Finally, we note that by applying the rescaling described in Section \ref{sec:rescale}, we may assume that the problem domain is the interval $[-1,1]$. After producing the solution, we can apply the inverse rescaling \eqref{eq:inverse_rescale} to return to the original domain $[a,b]$.

\subsection{The solution on $n$ nodes in a single subinterval} \label{sec:alg_nnodes}
In this section, we describe how to obtain a solution $\phi$ on $n$ Gaussian nodes of a subinterval $[x_{i},x_{i+1}]$. We will rescale the subinterval $[x_i,x_{i+1}]$ to be $[-1,1]$, as described in Section \ref{sec:rescale}; when the algorithm is completed, we will rescale back to the original domain $[x_i,x_{i+1}]$, as in equation \eqref{eq:inverse_rescale}. We explain how to produce the solution $\phi$, with specified boundary values, to the equation
\begin{math}
\L \phi = f
\end{math}
on $n$ Gaussian nodes $y_1,\dots, y_n$ in $[-1,1]$. We will also assume that we have divided the entire equation \eqref{eq:main_prob} by the leading coefficient $a_4(x)$, so that without any loss in generality we assume that the leading coefficient $a_4(x) = 1$.

As described in Section \ref{sec:fredholm}, we can write the solution $\phi$ as $\phi = \G_0 \sigma + \psi_\alpha$, where $\psi_\alpha$ is defined by formulas \eqref{eq:psi_alpha}--\eqref{eq:psi_r1}, and the function $\sigma$ satisfies the following second-kind integral equation:
\begin{align} \label{eq:int_eq000}
    \sigma(x) + \int_{-1}^{1} \sum_{j=0}^{3} a_j(x) G_j(x,t) \sigma(t) dt
      = f - \L \psi_\alpha \equiv f_{\alpha}.
\end{align}

We discretize the equation \eqref{eq:int_eq000} by sampling the coefficients $a_j(x)$, the right hand side $f_\alpha(x)$, and the functions $G_j(x,t)$ in each variable on the $n$ Gaussian nodes $y_1,\dots,y_n$.  We will produce the solution $\sigma$ evaluated on the same nodes. We will write the integrals in \eqref{eq:int_eq000} using the quadrature formula for Gaussian nodes and weights $w_k$ described in Section \ref{sec:legendre}. The discrete system of equations resulting from discretizing \eqref{eq:int_eq000} is given by:
\begin{align} \label{eq:syst000}
    \sigma(y_i) + \sum_{k=1}^{n} \left[ \sum_{j=0}^{3} 
        a_j(y_i) G_j(y_i,y_k) \right] \sigma(y_k) w_k = f_\alpha(y_i),
    \quad  1 \le i \le n.
\end{align}
We can compactly write this as a linear system $A \hat{\sigma} = \hat{f}_\alpha$, where we have defined $\hat{f}_\alpha = (f_\alpha(y_1), \dots, f_\alpha(y_n))^\top$, $\hat{\sigma} = (\sigma(y_1), \dots, \sigma(y_n))^\top$ and $A$ is the $n$-by-$n$ matrix with entries
\begin{align}
A(i,k) = \delta_{i,k} + \sum_{j=0}^{3} a_j(y_i) G_j(y_i,y_k) w_k.
\end{align}
We invert the $n$-by-$n$ matrix $A$ using the $QR$ factorization at a cost of $O(n^3)$ floating point operations \cite{dahlquist}; other algorithms for solving linear systems may be used as well.

The solution $\phi$ to \eqref{eq:main_prob}--\eqref{bdry4} and its first three derivatives on the nodes $y_1,\dots,y_n$ are now given by: 
\begin{align}
    \phi^{(j)}(y_i) = \sum_{k=1}^n G_j(y_i,y_k) \hat{\sigma}_k w_k 
                      + \psi_\alpha^{(j)}(y_i),
    \quad 0 \le j \le 3
\end{align}
and $\phi^{(4)}(y_i) = \hat{\sigma}_i$. We then perform the inverse rescaling in equation \eqref{eq:inverse_rescale} to complete the computation of the solution on the subinterval $[x_i,x_{i+1}]$. Since the cost of solving \eqref{eq:syst000} on each of the $m$ subintervals $[x_i,x_{i+1}]$ is $O(n^3)$, the total cost is $O(m  n^3)$.

\subsection{Matching the boundary values} \label{sec:alg_match}

In Section \ref{sec:alg_nnodes} we detailed how to obtain a solution to \eqref{eq:main_prob}--\eqref{bdry4} on $n$ Gaussian nodes in each subinterval $[x_i,x_{i+1}]$ with specified boundary values at $x_{i}$ and $x_{i+1}$. We will denote by $\tilde{\phi}_i$ the solution on $[x_i,x_{i+1}]$. In this section we explain how to adjust the boundary values of these functions so that the resulting functions $\phi_i$ and $\phi_{i+1}$ have matching values and three derivatives at the interface $x_{i+1}$, and so that $\phi_{1}$ satisfies \eqref{bdry1} and \eqref{bdry3} at $-1$, and $\phi_{m}$ satisfies \eqref{bdry2} and \eqref{bdry4} at $1$. We will assume that the $\tilde{\phi}_i$ were chosen to have zero boundary values at $x_{i}$ and $x_{i+1}$.

In addition to constructing $\tilde{\phi}_i$, on each subinterval $[x_i,x_{i+1}]$, we use the method of Section \ref{sec:alg_nnodes} to solve the equation $\L g = 0$ four times to obtain functions $g_{i,j}, 1 \le j \le 4$, with a single non-zero boundary value; concretely, $g_{i,1}(x_i) = 1, g_{i,2}(x_{i+1}) = 1, g_{i,3}^{\prime}(x_i)=1$, and $g_{i,4}^{\prime}(x_{i+1})=1$, and the other values and first derivatives at $x_i$ and $x_{i+1}$ are zero. We will find coefficients $\beta_{i,j}, 1 \le i \le m, 1 \le j \le 4$, such that the functions
\begin{align}
\label{adjust111}
\phi_i(x) = \tilde{\phi}_i(x) + \sum_{j=1}^{4} \beta_{i,j} g_{i,j}(x)
\end{align}
have matching boundary values, and the desired values at $\pm 1$ given by \eqref{bdry1}--\eqref{bdry4}. The derivatives of $\phi_i$ are found by adding the same linear combination of the derivatives of the $g_{i,j}$, namely:
\begin{align}
\label{adjust222}
\phi_i^{(k)}(x) = \tilde{\phi}_i^{(k)}(x) + \sum_{j=1}^{4} \beta_{i,j} g_{i,j}^{(k)}(x),
    \quad 1 \le k \le 4.
\end{align}
The function $\phi^{(k)}(x)$, $0 \le k \le 4$, is then defined by $\phi_i^{(k)}(x)$ when $x \in [x_i,x_{i+1}]$. In particular, $\sigma = \phi^{(4)}$ solves the integral equation \eqref{eq:int_eq000} on the entire interval $[-1,1]$.

\begin{figure}[h]
\centerline{
\includegraphics[width=0.5\linewidth]{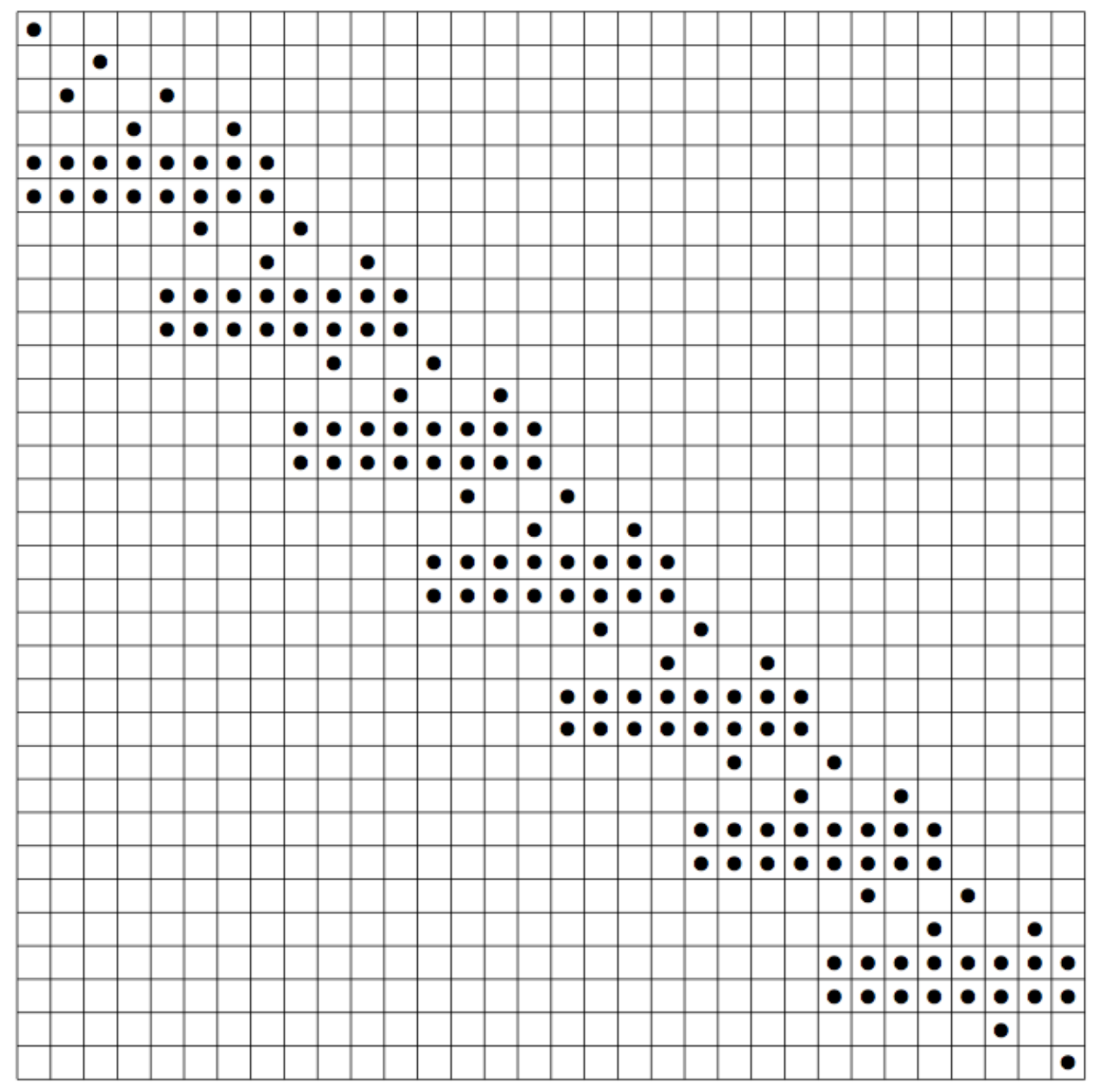}
}
\caption{The banded matrix used to match the boundary values, shown here for $m=8$ subintervals; non-zero values are marked with dots.}
\label{fig:band_matrix}

\end{figure}

The coefficients $\beta_{i,j}$ can be found as the solution to a banded linear system of
size $4m$-by-$4m$. For every pair of adjacent intervals $[x_{i}, x_{i+1}]$ and 
$[x_{i+1},x_{i+2}]$, we require that 
\begin{align}
\label{eq:phii_val}
\tilde{\phi}_{i}(x_{i+1}) + \beta_{i,2} 
     &=  \tilde{\phi}_{i+1}(x_{i+1}) + \beta_{i,1} \\ 
\label{eq:phii_der1}
\tilde{\phi}_{i}^{\prime}(x_{i+1}) + \beta_{i,4} 
     &=  \tilde{\phi}_{i+1}^{\prime}(x_{i+1}) + \beta_{i,3} \\ 
\label{eq:phii_der2}
\tilde{\phi}_{i}^{\prime \prime} (x_{i+1})+ \sum_{j=1}^{4} \beta_{i,j} g_{i,j}^{\prime 
        \prime}(x_{i+1})
&= \tilde{\phi}_{i+1}^{\prime \prime} (x_{i+1})+ \sum_{j=1}^{4} \beta_{i+1,j} 
        g_{i+1,j}^{\prime \prime}(x_{i+1}) \\
\label{eq:phii_der3}
\tilde{\phi}_{i}^{(3)} (x_{i+1})+ \sum_{j=1}^{4} \beta_{i,j} g_{i,j}^{(3)}(x_{i+1})
     &= \tilde{\phi}_{i+1}^{(3)} (x_{i+1})+ \sum_{j=1}^{4} \beta_{i+1,j} g_{i+1,j}^{(3)}
        (x_{i+1})
\end{align}
These conditions ensure that the functions $\phi_i$ have four matching derivatives at 
the interfaces $x_{i+1}, 1 \le i \le m-1$. To ensure the boundary conditions \eqref{bdry1}--\eqref{bdry4} on $[a,b]$, we also require the following equations:
\begin{align}
\label{eq:phi1_leftval}
\tilde{\phi}_{1}(-1) + \beta_{1,1} &= \alpha_{l,0} \\
\label{eq:phi1_rightval}
\tilde{\phi}_{1}^{\prime}(-1) + \beta_{1,3} &= \alpha_{l,1} \\
\label{eq:phim_leftder}
\tilde{\phi}_{m}(1) + \beta_{m,2} &= \alpha_{r,0} \\
\label{eq:phim_rightder}
\tilde{\phi}_{m}^{\prime}(1) + \beta_{m,4} &= \alpha_{r,1}.
\end{align}
Ordering the variables $\beta_{1,1},\beta_{1,2},\beta_{1,3},\beta_{1,4},\dots \dots,\beta_{m,1},\beta_{m,2},\beta_{m,3},\beta_{m,4}$, and ordering the equations as \eqref{eq:phii_val}--\eqref{eq:phii_der3} and by the intervals they involve, the linear system described by equations \eqref{eq:phii_val}--\eqref{eq:phim_rightder} is $9$-diagonal and so can be solved in $O(m)$ floating-point operations. In Figure \ref{fig:band_matrix}, we plot this matrix for $m=8$ subintervals. After solving for the coefficients $\beta_{i,j}$, the solutions $\phi_i$ on each subinterval can then be obtained by equations \eqref{adjust111} and \eqref{adjust222}, at an additional cost of $O(n  m)$.

\begin{rmk} \label{rmk:inaccurate}
Because the initial local solutions $\tilde{\phi}_i$ on each subinterval $[x_i,x_{i+1}]$ are obtained via the second-kind integral formulation described in Section \ref{sec:fredholm}, they can be computed with high accuracy. However, the method we have described in this section for matching their boundary values departs from this integral equation formulation. In practice, we have observed that the banded linear system defined by equations \eqref{eq:phii_der1}--\eqref{eq:phim_rightder} behaves like the matrices encountered in finite element or finite difference schemes, in that it has a condition number of size $O(m^4)$. To achieve greater accuracy, in Section \ref{sec:alg_corrections} we show how to apply deferred corrections to this problem.

\end{rmk}

\subsection{Applying deferred corrections} \label{sec:alg_corrections}
The solution we obtain from Sections \ref{sec:alg_nnodes} and \ref{sec:alg_match} will not be accurate to the precision permitted by the problem, as explained in Remark \ref{rmk:inaccurate}. To improve the accuracy, we apply the method of deferred corrections, as described in Section \ref{sec:corrections}. If $\hat{\sigma}$ is the vector of length $m  n$ with the estimated solution to the integral equation \eqref{eq:int_eq000} evaluated on the $m n$ nodes, then to apply deferred corrections we perform the following steps:
\begin{enumerate}

\item \label{apply_L}
Apply the operator $\L \G_0$ to $\hat{\sigma}$, evaluated on the $m n$ nodes.

\item \label{solve_new}
Solve the new system $\L \G_0 r = f_\alpha - \L \G_0 \hat{\sigma}$, obtaining an estimated residual $\hat{r}$.

\item \label{update}
Update the solution: $\hat{\sigma} \leftarrow \hat{\sigma} + \hat{r}$.
\end{enumerate}

Step \ref{solve_new} is performed by following Sections \ref{sec:alg_nnodes} and \ref{sec:alg_match} at a cost of $O(m  n^3)$ floating point operations, and Step \ref{update} requires $m n$ additions. We will now exhibit an algorithm for performing Step \ref{apply_L} with asymptotic cost $O(m  n^2)$.


%
Because $a_4 \equiv 1$, from equation \eqref{LG0} we have
\begin{align} 
(\L\G_0 \hat\sigma)(x) 
= \hat \sigma(x)
     + \sum_{j=0}^{3}a_j(x)  \int_{-1}^{1}  G_j(x,t) \hat \sigma(t) dt.
\end{align}
Therefore, we must integrate $\hat \sigma(x)$ against the functions $G_j(x,t)$. We now explain how to compute the integrals $\int_{-1}^{1} G_j(x,t) \hat \sigma(t) dt$ with $O(m  n^2 )$ floating-point operations. The function $G_j(x,t)$ can be written in the form
\begin{align} \label{eq:green_poly000}
    G_j(x,t) = 
    \begin{cases}
        \sum_{i=0}^{3-j} x^i p_i(t) , \text{ if } x < t \\
        \sum_{i=0}^{3-j} x^i q_i(t), \text{ if } t < x
    \end{cases}
\end{align}
where $p_i$ and $q_i$ are polynomials of degree $3-j$ that can be explicitly computed from
\eqref{eq:Green}. (Of course, $p_i$ and $q_i$ depend on $j$, though we suppress this for notational simplicity.) Using \eqref{eq:green_poly000}, we obtain:
\begin{align} \label{green_expand000}
    \int_{-1}^{1} G_j(x,t) \hat\sigma(t) dt
    = \sum_{i=0}^{3-j} x^{i} \bigg\{ \int_{-1}^{x} q_i(t) \hat\sigma(t) dt 
                                       + \int_{x}^{1} p_i(t) \hat\sigma(t) dt \bigg\} 
\end{align}
We must compute $\int_{-1}^{1} G_j(x,t) \hat\sigma(t) dt$ for each of the $m n $ values of $x$ on which we have discretized the problem, namely the $n$ Gaussian nodes on each of the $m$ subintervals $[x_{i},x_{i+1}]$. For each such $x$, let us suppose $x \in [x_{i^*},x_{i^* + 1}]$. We can then write:
\begin{align}
    \int_{-1}^{x} q_i(t) \hat \sigma(t) dt 
    = \sum_{k=1}^{i^*-1} \int_{x_{k}}^{x_{k+1}} q_i(t) \hat \sigma(t) dt 
                                    + \int_{x_{i^*}}^{x} q_i(t) \hat \sigma(t) dt 
\end{align}
and similarly
\begin{align}
    \int_{x}^{1} p_i(t) \hat \sigma(t) dt 
    = \sum_{k=i^*+1}^{m} \int_{x_{k}}^{x_{k+1}} p_i(t)\hat \sigma(t) dt 
                                   + \int_{x}^{x_{i^*+1}} p_i(t)\hat \sigma(t) dt .
\end{align}

Each of the $2m$ integrals $\int_{x_{k}}^{x_{k+1}} q_i(t) \hat\sigma(t) dt $ and $\int_{x_{k}}^{x_{k+1}} p_i(t) \hat\sigma(t) dt $ can be computed using an $n$-point Gaussian quadrature, at a total cost of $O(m n)$ floating-point operations. For each Gaussian node $x$ in the interval $[x_{i^*}, x_{i^*+1}]$, the integrals $\int_{x_{i^*}}^{x} q_i(t) \hat\sigma(t) dt$ and $\int_{x}^{x_{i^*+1}} p_i(t) \hat\sigma(t) dt$ can also be computed using an $n$-point Gaussian quadrature, at a total cost of $O(m n^2)$ floating-point operations.

From \eqref{green_expand000}, the integrals $\int_{-1}^1 G_j(x,t) \hat{\sigma}(t) dt$, and hence $(\L \G_0 \hat{\sigma})(x)$, can be computed at cost $O(m  n^2)$. Consequently, each iteration of deferred corrections requires $O(mn^2)$ operations.  As explained in Section \ref{sec:corrections}, the number of iterations grows like $O(\log(1/\ep))$, where $\ep$ is the machine precision. This brings the asymptotic running time of the entire algorithm to $O(m n^3 \log(1/\ep))$.

\section{Numerical results} \label{sec:experiments}
The algorithm of Section \ref{sec:algorithm} has been implemented in Fortran. In this section, we use the algorithm to obtain the numerical solution $\phi$ to several specific problems of the form \eqref{eq:main_prob}--\eqref{bdry4}. All experiments were performed on a Dell XPS-L521X laptop computer with an Intel Core i7-3632QM CPU at 2.20GHz with 15.6 GB RAM, running the 64-bit Ubuntu 14.04 LTS operating system. In our experiments, unless otherwise specified the code was compiled with the gfortran compiler using extended (16-bit) precision; this is called with the compilation flag \texttt{-freal-8-real-16}. Since most users are likely to use the method with double precision, our use of the extended precision environment leads us to report somewhat pessimistic CPU times. We have made no attempt to optimize the choice of the parameter $n$ to achieve maximal performance.


Six experiments are presented below, and their results are contained in Tables \ref{errors:sin150}--\ref{residuals:bdlayer}. In each experiment, we solve a specific boundary value problem on a specified interval $[a,b]$, and measure the relative error of the solution and its first four derivatives on a grid of $N=10,000$ equispaced points $x_1,\dots,x_N$ on $[a,b]$, including the endpoints $a$ and $b$. If $\hat{\phi}^{(j)}$ is the estimated $j^{th}$ derivative, then the relative error is defined as follows:
\begin{align}
    R(\phi^{(j)}) = \sqrt{ 
           \frac{\sum_{i=1}^{N} | \hat{\phi}^{(j)}(x_i) - \phi^{(j)}(x_i) |^2}
                {\sum_{i=1}^{N} |\phi^{(j)}(x_i) |^2}  }
\end{align}


For each experiment, we display two tables. In each table, the first column contains the number $m$ of subintervals of $[a,b]$ used for that experiment. In the first table, the first five columns after $m$ display the errors of the solution and its first four derivatives; the last column displays the running time, in seconds. In the second table, the columns after $m$ display the relative size of the residual norm after each iteration of the algorithm, to demonstrate the convergence of deferred corrections. In our implementation of the algorithm, we run several iterations after convergence, whose errors are approximately the same; to save space in the tables, these are not displayed.

\begin{rmk}
The relative errors $R(\phi^{(j)})$ decrease by approximately a fixed number of digits whenever $m$ is doubled. This is the behavior expected from the error of the
Gaussian quadrature, as described in Section \ref{sec:legendre}, equation  \eqref{leg_err}.

\end{rmk}

\begin{rmk}
For each example, the running times scale approximately linearly with $m$. This behavior is expected from the $O(m)$ asymptotic running time of the algorithm.
\end{rmk}

\begin{rmk}
In the examples, the residual size typically decreases by an approximately constant number of digits after approximately two iterations of deferred corrections. This is likely because two linear systems are solved each step -- the local linear systems on each subinterval, and the global linear system that matches the local solutions' boundary values. Up to taking two steps instead of one, this rate of decrease in the residual size is what we expect from Section \ref{sec:corrections}. In particular, very few iterations of the algorithm are required until convergence is achieved.
\end{rmk}

\subsection{The function $\phi(x) = \sin(150x)$}
We solve the following equation:
\begin{align}
    \sum_{j=0}^4 (1 + x^{4-j}) \phi^{(j)}(x) 
        = \sum_{j=0}^4 (1 + x^{4-j}) \frac{d^j}{dx^j} \sin(150x)
\end{align}
on the interval $[0,2\pi]$, subject to the boundary conditions $\alpha_{l,0}=0$, $\alpha_{r,0}=0$, $\alpha_{l,1}=150$ and $\alpha_{r,1}=150$. The solution is apparently $\phi(x) = \sin(150x)$. 
Because of the high frequency, we do not expect to resolve the function to full machine precision (regardless of the choice of algorithm); however, our algorithm successfully recovers the solution to within the error permitted by the problem's condition number. The algorithm was run using $n=15$ Gaussian nodes per subinterval. The results are displayed in Tables \ref{errors:sin150} and \ref{residuals:sin150}.

\subsection{Beam with fixed ends} \label{sec:fixed_ends}
In this example, we consider the problem of determining the shape of a beam with non-uniform stiffness subjected to an external force. We parametrize the $x$-axis by the interval $[0,1]$. We take the stiffness of the beam to be the function $c(x)$ given by:
\begin{align}
    c(x) = (x-1/2)^2 + 1
\end{align}
and the external force to be the function $f(x)$ given by:
\begin{align}
    f(x) = \sin(2 \pi x) + 1.
\end{align}

\begin{figure}[h]
\centerline{
\includegraphics[width=.8\linewidth]{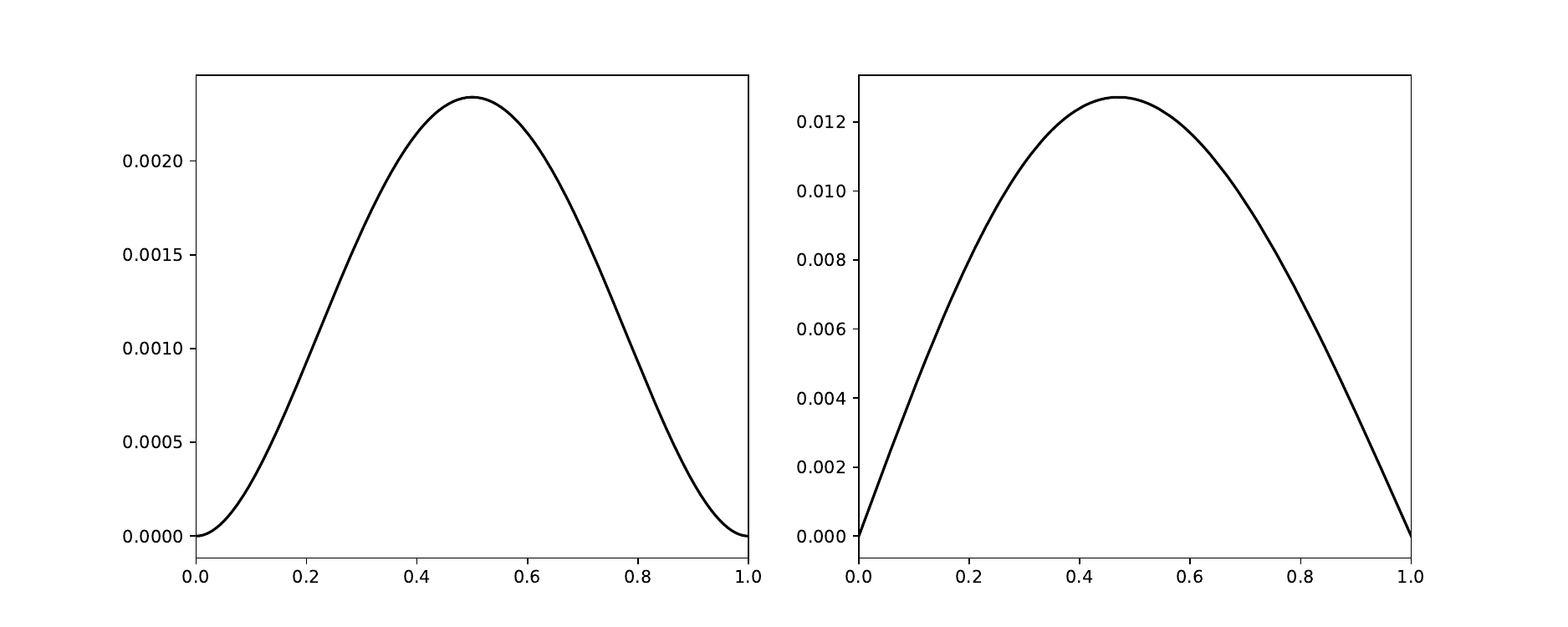}
}
\caption{Left: the beam with fixed ends. Right: the beam with simply-supported ends.}
\label{fig:beams}

\end{figure}

If $\phi(x)$ is the shape of the beam, then $\phi$ satisfies the following differential equation \cite{gere-mechanics-1984}:
\begin{align}
    \frac{d^2}{dx^2} \left( c(x) \frac{d^2 \phi}{dx^2} \right) = f(x)
\end{align}
or equivalently
\begin{align}  \label{spelled_out}
    ((x-1/2)^2 + 1) \frac{d^4 \phi}{dx^4}(x) + 4(x-1/2)\frac{d^3 \phi}{dx^3}(x)
           + 2 \frac{d^2 \phi}{dx^2}(x) = \sin(2 \pi x) + 1.
\end{align}
We impose fixed endpoints on the beam, meaning that $\alpha_{l,0}=0$, $\alpha_{r,0}=0$, $\alpha_{l,1}=0$ and $\alpha_{r,1}=0$. The algorithm was run using $n=10$ Gaussian nodes per subinterval. The results are displayed in Tables \ref{errors:fixed} and \ref{residuals:fixed}.  The solution is plotted on the left side of Figure \ref{fig:beams}.

\subsection{Beam with simply-supported ends}
\label{sec-simply}
In this example, we consider the same beam-bending problem \eqref{spelled_out} as in Section \ref{sec:fixed_ends}, except we require that the ends of the beam are simply-supported, or $\phi(0) = \phi(1) = 0$ and $\phi^{\prime\prime}(0) = \phi^{\prime\prime}(1) = 0$. Because the boundary vales involve the second derivatives, this problem does not immediately fit into the form of \eqref{eq:main_prob}--\eqref{bdry4}, and we include it as an example for that reason.

To introduce the correct boundary conditions, we let $\L$ denote the differential operator on the left side of \eqref{spelled_out}, and we solve the equation $\L \phi = f$ with boundary conditions $\alpha_{l,0} = \alpha_{r,0} = \alpha_{l,1} = \alpha_{r,1} = 0$; we will call this function $\tilde{\phi}$. We also produce four linearly independent solutions to the equation $\L \phi = 0$ by solving it with linearly independent boundary conditions; we will call these solutions $\phi_{1}$, $\phi_2$, $\phi_3$ and $\phi_4$.

We will find coefficients $b_i$, $1 \le i \le 4$, so that the function
\begin{align} \label{phi_plus_bphi}
    \phi(x) = \tilde \phi(x) + \sum_{i=1}^4 b_i \phi_i(x)
\end{align}
satisfies the sought-after boundary conditions $\phi(0) = \phi(1) = 0$ and $\phi^{\prime\prime}(0) = \phi^{\prime\prime}(1) = 0$. The $b_i$ satisfy the linear equations:
\begin{align}
    0 &= \sum_{i=1}^4 b_i \phi_i(0) \\
    0 &= \sum_{i=1}^4 b_i \phi_i(1) \\
    0 &= \tilde{\phi}^{\prime\prime}(0) + \sum_{i=1}^4 b_i \phi_i^{\prime\prime}(0) \\
    0 &= \tilde{\phi}^{\prime\prime}(1) + \sum_{i=1}^4 b_i \phi_i^{\prime\prime}(1) . 
\end{align}

This system is non-singular because the $\phi_i$ are linearly independent. This permits us to  solve for the $b_i$, and then define the solution $\phi$ by \eqref{phi_plus_bphi}. The algorithm was run using $n=10$ Gaussian nodes per subinterval. The results are shown in Tables \ref{errors:simply} and \ref{residuals:simply}. The solution is plotted on the right side of Figure \ref{fig:beams}.

\begin{figure}[h]
\centerline{
\includegraphics[width=0.5\linewidth]{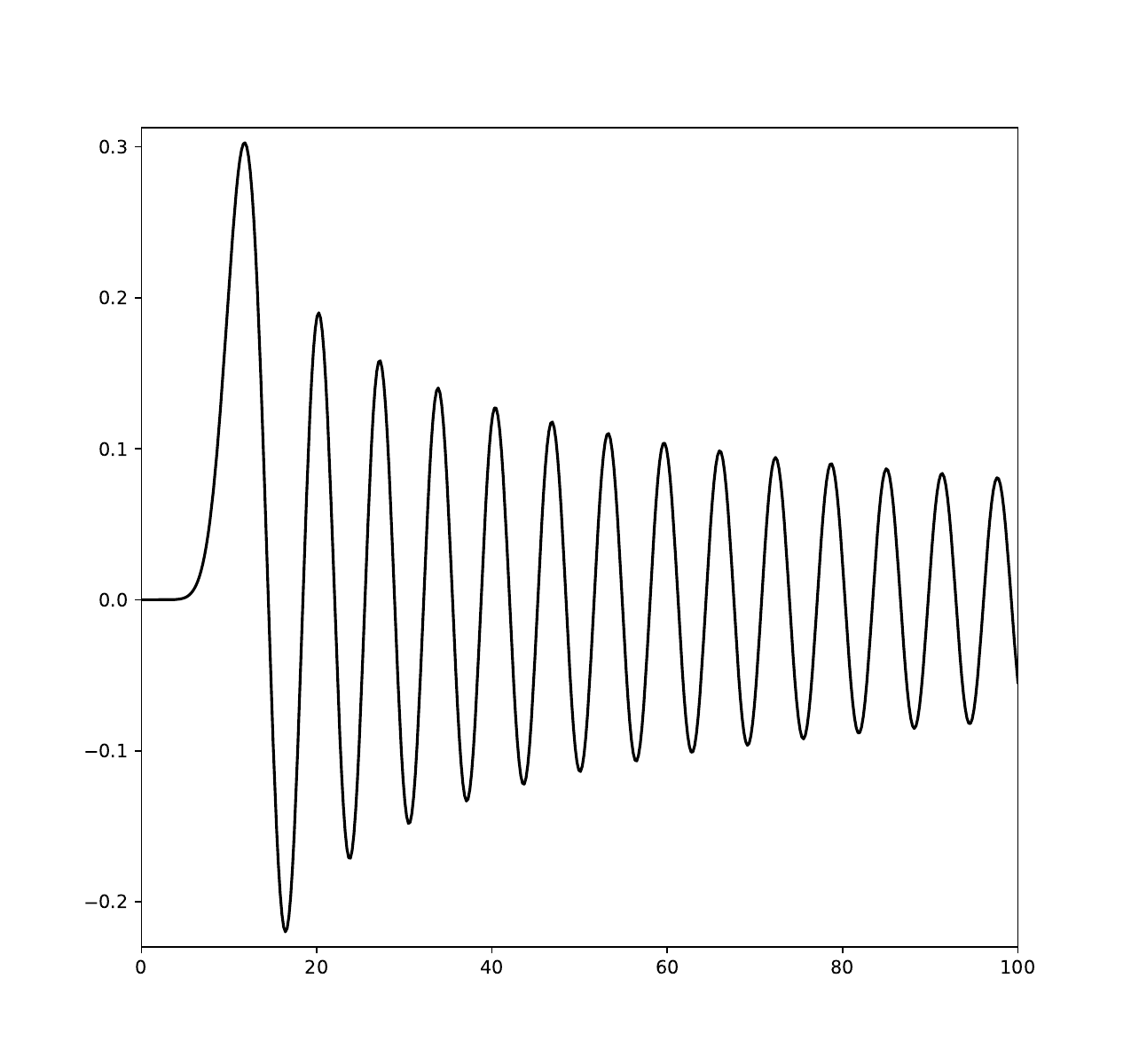}
}
\caption{The Bessel function $J_{10}$ on the interval $[0,100]$.}
\label{fig:bessel}

\end{figure}

\subsection{The Bessel function $J_{10}(x)$}

In this example, we solve for the Bessel function $J_{10}(x)$ on the interval $[0,100]$, which is plotted in Figure \ref{fig:bessel}. $J_{10}$ satisfies the second-order equation:
\begin{align}
    x^2 \frac{d^2}{dx^2} \phi(x) + x \frac{d}{dx} \phi(x) + (x^2-100) \phi(x) = 0.
\end{align}

We differentiate this equation twice to arrive at the fourth-order equation:
\begin{align}
\label{eq:bessel}
    \left(x^2 \frac{d^4}{dx^4} + 5x \frac{d^3}{dx^3}  + (x^2-96) \frac{d^2}{dx^2} 
        + 4x  \frac{d}{dx} + 2 \right) \phi(x)= 0.
\end{align}
We solve this equation subject to the boundary values matching the known values of $J_{10}(x)$ and $J_{10}^\prime(x)$ at 0 and 100. Because of the singularity at $x=0$, we represent the 0 endpoint by the square root of machine zero. We ran the algorithm with $n=20$ nodes per subinterval. We include this example to demonstrate the algorithm's performance when the leading coefficients has a root, as does the leading coefficient $x^2$ at the left endpoint $x=0$. Because of this singularity, the solution cannot be obtained to full machine precision even on the $m  n$ Gaussian nodes themselves. Nevertheless, our algorithm obtains the solution up to this necessary loss of accuracy. The results are shown in Tables \ref{errors:bessel} and \ref{residuals:bessel}.

\subsection{The function $\phi(x) = \exp(\sin(2x))$}
In this example, we consider the function $\phi(x) = \exp(\sin(2x))$ on the interval $[0,2\pi]$, which was examined in the paper \cite{driscoll2010automatic}. This function satisfies the differential equation
\begin{align}
\phi^{(4)}(x) - 2\cos(2x) \phi^{(3)}(x) 
         + & [48\cos^2(2x) (1 + \sin(2x)) 
    \nonumber \\
&        - 16\sin(2x)(1+3\sin(2x))] \phi(x) = 0,
\end{align}
with boundary conditions $\phi(0)=1$, $\phi^\prime(0)=2$, $\phi^\prime(2 \pi)=2$ and $\phi^{\prime\prime}(2 \pi)=4$. Since the right boundary condition involves the second derivative, we use the method described in Section \ref{sec-simply} to ensure the correct boundary values. We run the algorithm using $n=12$ Gaussian nodes per subinterval. The results are displayed in Tables \ref{errors:driscoll} and \ref{residuals:driscoll}.

We note that the method in \cite{driscoll2010automatic} is reported to take 0.37 seconds to achieve relative error $0.46 \times 10^{-12}$; by contrast, our method takes between 0.15 and 0.23 seconds to achieve comparable accuracy, running in an extended precision environment. Furthermore, in a separate test we run our algorithm in the double precision environment with $m=312$ subintervals and $n=7$ Gaussian nodes per subinterval; we achieve a relative error of $0.44 \times 10^{-12}$ in only 0.056 seconds. We caution against over-interpreting this comparison, however, as the method from \cite{driscoll2010automatic} was implemented in MATLAB, and adaptively finds the discretization; any comparison of the timings must be understood in this context.

\subsection{Function with a boundary layer}

We consider the example, adapted from \cite{greengard-two-point}, of the function on $[-1,1]$ given by the formula
\begin{align}
\phi(x) = A e^{x/\epsilon} + B,
\end{align}
where $A = 1/(e^{1/\epsilon} - e^{-1/\epsilon})$ and $B = (2e^{-1/\epsilon} - e^{1/\epsilon})/(e^{-1/\epsilon} - e^{1/\epsilon})$. This function satisfies the equation
\begin{align}
\epsilon \phi^{(4)}(x) - \phi^{(3)}(x)  = 0
\end{align}
with boundary conditions $\phi(-1) = 1$, $\phi(1) = 2$, $\phi'(-1) = A e^{-1/\epsilon}/\epsilon$, and $\phi'(1) = A e^{1/\epsilon}/\epsilon$.

When $\epsilon > 0$ is small, $\phi(x)$ has a sharp boundary layer at $1$, and most standard techniques are likely to fail for this reason. In our experiments, we take $\epsilon = 10^{-4}$. Because of the boundary layer, we use a grid whose interval lengths become geometrically smaller as they approach the right end point. In more detail, we construct $m'$ subintervals $[x_i',x_{i+1}']$ defined by 
\begin{align}
& x_1' = -1;
\nonumber \\
& x_{i+1}' = x_{i}' + 2^{-i+1}, \quad 1 \le i \le m';
\nonumber \\
& x_{m' + 1}' = 1.
\end{align}
We then equally subdivide each of these intervals by $10$, for a total of $m = 10 m'$ subintervals. We use $n=28$ Gaussian nodes per subinterval. As with the Bessel function, the ill-conditioning of the problem prevents the function from being fully resolved, even on the $mn$ Gaussian nodes themselves. Nevertheless, the algorithm recovers the solution up to this necessary loss of accuracy. The results are displayed in Tables \ref{errors:bdlayer} and \ref{residuals:bdlayer}.

\section*{Acknowledgements}

William Leeb acknowledges support from a postdoctoral fellowship from the Simons Collaboration on Algorithms and Geometry, NSF award IIS-1837992, and BSF award 2018230. Vladimir Rokhlin acknowledges support from the Office of Naval Research under the grant N00014-16-1-2123 and the Air Force Office of Scientific Research under the grant FA9550-16-1-0175. The authors thank the referees for their helpful comments on the manuscript.


\bibliographystyle{plain}
\bibliography{./mybib}

\newpage

%
%
%
%
 \begin{table}
 \caption{Relative errors and running times for                            $\sin(150x)$}
 \label{errors:sin150} 
 \small 
 \center 
 \begin{tabular}{|l| c  | c | c | c | c | c |}  
 \hline
 $m$ & $R(\phi)$ & $R(\phi^\prime)$ &                     $R(\phi^{\prime\prime})$ & $R(\phi^{(3)})$  &                     $R(\phi^{(4)})$   &  Time (s)\\  
 \hline
   16  &      0.10E+01  &      0.10E+01  &      0.10E+01  &      0.10E+01  &      0.13E+01  
        & 0.17E+00\\
   32  &      0.43E+00  &      0.36E+00  &      0.36E+00  &      0.37E+00  &      0.34E+00  
        & 0.31E+00 \\
   64  &      0.13E-03  &      0.13E-03  &      0.14E-03  &      0.14E-03  &      0.13E-03  
        & 0.57E+00  \\
  128  &      0.74E-08  &      0.74E-08  &      0.74E-08  &      0.74E-08  &      0.74E-08  
        &  0.11E+01 \\
  256  &      0.26E-12  &      0.26E-12  &      0.26E-12  &      0.26E-12  &      0.26E-12  
        &  0.22E+01 \\
  512  &      0.84E-17  &      0.82E-17  &      0.84E-17  &      0.82E-17  &      0.84E-17  
        &   0.43E+01 \\
 1024  &      0.26E-21  &      0.25E-21  &      0.26E-21  &      0.25E-21  &      0.26E-21  
        &   0.85E+01\\
 2048  &      0.68E-25  &      0.78E-26  &      0.78E-26  &      0.78E-26  &      0.78E-26  
        & 0.16E+02\\
 4096  &      0.31E-25  &      0.24E-27  &      0.15E-29  &      0.25E-30  &      0.72E-30  
        & 0.32E+02\\
 \hline
 \end{tabular}
 \end{table} 
 \begin{table} 
 \caption{Relative residual sizes for                    $\sin(150x)$}
 \label{residuals:sin150}.
 \small 
 \center 
 \begin{tabular}{|l |c |c |c |c |c |c |c | } 
 \hline
& \multicolumn{    7             }{c|}{Iteration number}\\   \cline{2-    8 }
 $m$ &   1 &   2 &   3 &   4 &   5 &   6 &   7 \\
 \hline
   16  &      0.91E-02  &      0.16E-05  &      0.16E-20  &      0.17E-22  &      0.30E-33  &      0.52E-33  &      0.31E-33  \\  
   32 &     0.30E-01 &     0.72E-06 &     0.23E-24 &     0.34E-27 &     0.21E-33 &     0.24E-33 &    \\ 
   64 &     0.13E-04 &     0.32E-09 &     0.25E-32 &     0.24E-33 &     0.13E-33 &    &    \\ 
  128 &     0.23E-08 &     0.93E-14 &     0.16E-33 &     0.17E-33 &     0.21E-33 &    &    \\ 
  256 &     0.13E-13 &     0.11E-17 &     0.62E-33 &     0.21E-33 &     0.14E-33 &    &    \\ 
  512 &     0.12E-17 &     0.26E-22 &     0.44E-33 &     0.14E-33 &     0.19E-33 &    &    \\ 
 1024 &     0.95E-22 &     0.44E-27 &     0.24E-33 &     0.18E-33 &     0.19E-33 &    &    \\ 
 2048 &     0.62E-26 &     0.68E-32 &     0.24E-33 &     0.14E-33 &    &    &    \\ 
 4096 &     0.40E-30 &     0.17E-33 &     0.16E-33 &    &    &    &    \\ 
 \hline
 \end{tabular}
 \end{table}

%
%
%
%
 \begin{table} 
 \caption{Relative errors and running times for                            beam with fixed ends}
 \label{errors:fixed}
 \small 
 \center 
 \begin{tabular}{|l| c  | c | c | c | c | c |}  
 \hline
 $m$ & $R(\phi)$ & $R(\phi^\prime)$ &                     $R(\phi^{\prime\prime})$ & $R(\phi^{(3)})$  &                     $R(\phi^{(4)})$   &  Time (s)\\  
 \hline
    2  &      0.27E-07  &      0.74E-07  &      0.66E-07  &      0.18E-06  &      0.13E-06  
        &      0.20E-01\\
    4  &      0.30E-10  &      0.51E-10  &      0.14E-09  &      0.51E-10  &      0.68E-09  
         &      0.24E-01\\
    8  &      0.27E-13  &      0.71E-13  &      0.11E-12  &      0.15E-12  &      0.49E-12  
         &      0.40E-01\\
   16  &      0.26E-16  &      0.70E-16  &      0.11E-15  &      0.15E-15  &      0.49E-15  
         &      0.68E-01\\
   32  &      0.26E-19  &      0.68E-19  &      0.10E-18  &      0.15E-18  &      0.48E-18  
           &      0.11E+00\\
   64  &      0.25E-22  &      0.66E-22  &      0.10E-21  &      0.15E-21  &      0.47E-21  
         &      0.22E+00\\
  128  &      0.24E-25  &      0.65E-25  &      0.10E-24  &      0.14E-24  &      0.45E-24  
          &      0.41E+00\\
  256  &      0.24E-28  &      0.63E-28  &      0.97E-28  &      0.14E-27  &      0.44E-27  
           &      0.82E+00\\
  512  &      0.70E-31  &      0.11E-30  &      0.14E-30  &      0.14E-30  &      0.43E-30 
            &      0.16E+01 \\
 1024  &      0.66E-31  &      0.92E-31  &      0.11E-30  &      0.28E-31  &      0.13E-31 
             &      0.32E+01 \\
 \hline
 \end{tabular}
 \end{table}

 \begin{table} 
 \caption{Relative residual sizes for                    beam with fixed ends}
 \label{residuals:fixed}
 \small 
 \center 
 \begin{tabular}{|l |c |c |c |c |c |c |c | } 
 \hline
& \multicolumn{    7             }{c|}{Iteration number}\\   \cline{2-    8 }
 $m$ &   1 &   2 &   3 &   4 &   5 &   6 &   7 \\
 \hline
    2  &      0.51E-07  &      0.41E-07  &      0.68E-15  &      0.19E-15  &      0.53E-23  &      0.88E-24  &      0.34E-31  \\  
    4 &     0.83E-10 &     0.21E-10 &     0.53E-21 &     0.60E-22 &     0.25E-32 &    &    \\ 
    8 &     0.12E-12 &     0.51E-14 &     0.22E-27 &     0.54E-29 &     0.11E-33 &    &    \\ 
   16 &     0.12E-15 &     0.15E-17 &     0.29E-33 &     0.91E-34 &     0.61E-34 &    &    \\ 
   32 &     0.11E-18 &     0.42E-20 &     0.12E-33 &    &    &    &    \\ 
   64 &     0.11E-21 &     0.53E-23 &     0.12E-33 &    &    &    &    \\ 
  128 &     0.11E-24 &     0.57E-26 &     0.12E-33 &    &    &    &    \\ 
  256 &     0.67E-27 &     0.15E-28 &     0.22E-33 &    &    &    &    \\ 
  512 &     0.88E-26 &     0.35E-28 &     0.27E-33 &    &    &    &    \\ 
 1024 &     0.12E-25 &     0.56E-28 &     0.30E-33 &    &    &    &    \\ 
 \hline
 \end{tabular}
 \end{table}

%
%
%
%
 \begin{table} 
 \caption{Relative errors and running times for                            beam with simply-supported ends}
 \label{errors:simply}
 \small 
 \center 
 \begin{tabular}{|l| c  | c | c | c | c | c |}  
 \hline
 $m$ & $R(\phi)$ & $R(\phi^\prime)$ &                     $R(\phi^{\prime\prime})$ & $R(\phi^{(3)})$  &                     $R(\phi^{(4)})$  & Time (s)\\  
 \hline
    2  &      0.29E-07  &      0.32E-07  &      0.49E-07  &      0.17E-06  &      0.12E-06  
         &      0.28E-01\\
    4  &      0.28E-10  &      0.43E-10  &      0.11E-09  &      0.10E-09  &      0.61E-09  
          &      0.40E-01\\
    8  &      0.28E-13  &      0.44E-13  &      0.10E-12  &      0.17E-12  &      0.44E-12  
           &      0.72E-01\\
   16  &      0.26E-16  &      0.39E-16  &      0.92E-16  &      0.16E-15  &      0.44E-15 
           &      0.13E+00 \\
   32  &      0.24E-19  &      0.35E-19  &      0.84E-19  &      0.15E-18  &      0.43E-18 
           &      0.19E+00 \\
   64  &      0.23E-22  &      0.33E-22  &      0.79E-22  &      0.15E-21  &      0.42E-21 
           &      0.38E+00  \\
  128  &      0.22E-25  &      0.32E-25  &      0.76E-25  &      0.14E-24  &      0.41E-24 
           &      0.80E+00 \\
  256  &      0.21E-28  &      0.31E-28  &      0.73E-28  &      0.14E-27  &      0.40E-27
           &      0.14E+01  \\
  512  &      0.16E-30  &      0.11E-30  &      0.14E-30  &      0.14E-30  &      0.39E-30 
          &      0.28E+01 \\
 1024  &      0.16E-30  &      0.11E-30  &      0.10E-30  &      0.46E-31  &      0.18E-31  
            &      0.54E+01\\
 \hline
 \end{tabular}
 \end{table}
 \begin{table} 
 \caption{Relative residual sizes for                    beam with simply-supported ends}
 \label{residuals:simply}
 \small 
 \center 
 \begin{tabular}{|l |c |c |c |c |c |c |c | } 
 \hline
& \multicolumn{    7             }{c|}{Iteration number}\\   \cline{2-    8 }
 $m$ &   1 &   2 &   3 &   4 &   5 &   6 &   7 \\
 \hline
    2  &      0.51E-07  &      0.41E-07  &      0.68E-15  &      0.19E-15  &      0.53E-23  &      0.88E-24  &      0.34E-31  \\  
    4 &     0.83E-10 &     0.21E-10 &     0.53E-21 &     0.60E-22 &     0.25E-32 &    &    \\ 
    8 &     0.12E-12 &     0.51E-14 &     0.22E-27 &     0.54E-29 &     0.11E-33 &    &    \\ 
   16 &     0.12E-15 &     0.15E-17 &     0.29E-33 &     0.91E-34 &     0.61E-34 &    &    \\ 
   32 &     0.11E-18 &     0.42E-20 &     0.12E-33 &    &    &    &    \\ 
   64 &     0.11E-21 &     0.53E-23 &     0.12E-33 &    &    &    &    \\ 
  128 &     0.11E-24 &     0.57E-26 &     0.12E-33 &    &    &    &    \\ 
  256 &     0.67E-27 &     0.15E-28 &     0.22E-33 &    &    &    &    \\ 
  512 &     0.88E-26 &     0.35E-28 &     0.27E-33 &    &    &    &    \\ 
 1024 &     0.12E-25 &     0.56E-28 &     0.30E-33 &    &    &    &    \\ 
 \hline
 \end{tabular}
 \end{table}

%
%
%
%
 \begin{table} 
 \caption{Relative errors and running times for                            $J_{10}(x)$}
 \label{errors:bessel}
 \small 
 \center 
 \begin{tabular}{|l| c  | c | c | c | c | c |}  
 \hline
 $m$ & $R(\phi)$ & $R(\phi^\prime)$ &                     $R(\phi^{\prime\prime})$ & $R(\phi^{(3)})$  &                     $R(\phi^{(4)})$ & Time (s)\\  
 \hline
   16  &      0.21E-14  &      0.12E-14  &      0.28E-14  &      0.20E-13  &      0.82E-12 
          &      0.43E+00 \\
   32  &      0.10E-20  &      0.27E-20  &      0.15E-20  &      0.24E-18  &      0.19E-16  
         &      0.66E+00\\
   64  &      0.17E-26  &      0.81E-25  &      0.55E-24  &      0.42E-21  &      0.67E-19
         &      0.12E+01  \\
  128  &      0.26E-27  &      0.11E-27  &      0.15E-26  &      0.27E-23  &      0.84E-21
         &      0.24E+01  \\
  256  &      0.67E-27  &      0.30E-28  &      0.24E-28  &      0.91E-25  &      0.57E-22 
         &      0.46E+01 \\
  512  &      0.54E-27  &      0.25E-28  &      0.21E-28  &      0.53E-26  &      0.66E-23  
          &      0.92E+01\\
 \hline
 \end{tabular}
 \end{table} 
 \begin{table} 
 \caption{Relative residual sizes for                    $J_{10}(x)$}
\label{residuals:bessel}
 \small 
 \center
 \begin{tabular}{|l |c |c |c |c |c |c |c | } 
 \hline
& \multicolumn{    7             }{c|}{Iteration number}\\   \cline{2-    8 }
 $m$ &   1 &   2 &   3 &   4 &   5 &   6 &   7 \\
 \hline
   16  &      0.34E-10  &      0.12E-11  &      0.23E-12  &      0.42E-13  &      0.79E-14  &      0.15E-14  &      0.27E-15  \\  
   32 &     0.27E-16 &     0.76E-18 &     0.60E-19 &     0.64E-20 &     0.69E-21 &     0.74E-22 &     0.79E-23  \\  
   64 &     0.14E-19 &     0.15E-20 &     0.13E-21 &     0.12E-22 &     0.11E-23 &     0.93E-25 &     0.78E-26  \\  
  128 &     0.25E-22 &     0.37E-23 &     0.30E-24 &     0.24E-25 &     0.21E-26 &     0.60E-27 &    \\  
  256 &     0.21E-24 &     0.18E-25 &     0.13E-26 &     0.56E-27 &     0.41E-27 &  &     \\  
  512 &     0.48E-25 &     0.21E-27  &   &   &   &   &   \\  
 \hline
 \end{tabular}

\bigskip

 \begin{tabular}{|l |c |c |c |c |c |c |c | } 
 \hline
& \multicolumn{    7             }{c|}{Iteration number}\\   \cline{2-    8 }
 $m$ &   8 &   9 &   10 &   11 &   12 &   13 &   14 \\
 \hline
   16 \,   &      0.50E-16  &      0.92E-17  &      0.17E-17  &      0.32E-18  &      0.59E-19  &      0.11E-19  &      0.20E-20  \\  
   32  & 0.84E-24 &     0.91E-25 &     0.97E-26 &     0.15E-26 &     0.65E-28  &  & \\  
   64  &     0.36E-27 &     0.58E-28 &  & &  &  & \\  
 \hline
 \end{tabular}

\bigskip

 \begin{tabular}{|l |c |c |c |c |c |c | c |} 
 \hline
& \multicolumn{    7             }{c|}{Iteration number}\\   \cline{2-    8 }
 $m$ &   15 &   16 &   17 &   18 &   19 &   20  & 21 \\
 \hline
   16 \,  &     0.37E-21  &      0.69E-22  &      0.13E-22  &      0.24E-23  &      0.43E-24  &      0.90E-25  &      0.62E-26 \\  
 \hline
 \end{tabular}

 \end{table}

%
%
%
%
 \begin{table} 
 \caption{Relative errors and running times for                            $\exp(\sin(2x))$}
 \label{errors:driscoll}.
 \small 
 \center 
 \begin{tabular}{|l| c  | c | c | c | c | c |}  
 \hline
 $m$ & $R(\phi)$ & $R(\phi^\prime)$ &                     $R(\phi^{\prime\prime})$ & $R(\phi^{(3)})$  &                     $R(\phi^{(4)})$ & Time (s)\\  
 \hline
   16  &      0.10E-10  &      0.14E-10  &      0.28E-10  &      0.89E-10  &      0.27E-09  
           &      0.15E+00\\
   32  &      0.63E-14  &      0.78E-14  &      0.97E-14  &      0.26E-13  &      0.59E-13
          &      0.23E+00  \\
   64  &      0.18E-17  &      0.23E-17  &      0.27E-17  &      0.63E-17  &      0.17E-16
          &      0.43E+00  \\
  128  &      0.45E-21  &      0.56E-21  &      0.66E-21  &      0.15E-20  &      0.41E-20
          &      0.80E+00  \\
  256  &      0.11E-24  &      0.13E-24  &      0.16E-24  &      0.38E-24  &      0.10E-23
          &      0.14E+01  \\
  512  &      0.29E-28  &      0.35E-28  &      0.41E-28  &      0.92E-28  &      0.25E-27
          &      0.28E+01  \\
 1024  &      0.95E-29  &      0.10E-28  &      0.95E-29  &      0.95E-29  &      0.96E-29
          &      0.56E+01  \\
 \hline
 \end{tabular}
 \end{table} 
 \begin{table} 
 \caption{Relative residual sizes for                    $\exp(\sin(2x))$}
 \label{residuals:driscoll}.
 \small 
 \center 
 \begin{tabular}{|l |c |c |c |c |c |c | } 
 \hline
& \multicolumn{    6             }{c|}{Iteration number}\\   \cline{2-    7 }
 $m$ &   1 &   2 &   3 &   4 &   5 &   6 \\
 \hline
   16  &      0.51E-06  &      0.28E-08  &      0.26E-16  &      0.17E-18  &      0.13E-26  &      0.13E-28  \\  
   32 &     0.18E-10 &     0.79E-12 &     0.29E-24 &     0.15E-25 &     0.14E-28 &    \\ 
   64 &     0.54E-14 &     0.17E-15 &     0.13E-28 &     0.24E-28 &     0.16E-28 &    \\ 
  128 &     0.13E-17 &     0.35E-19 &     0.11E-28 &     0.11E-28 &    &    \\ 
  256 &     0.32E-21 &     0.75E-23 &     0.52E-29 &    &    &    \\ 
  512 &     0.31E-22 &     0.50E-26 &     0.22E-28 &    &    &    \\ 
 1024 &     0.21E-21 &     0.30E-25 &     0.42E-28 &    &    &    \\ 
 \hline
 \end{tabular}
 \end{table} 

%
%
%
%
 \begin{table} 
 \caption{Relative errors and running times for                            boundary layer}
 \label{errors:bdlayer}.
 \small 
 \center 
 \begin{tabular}{|l| c  | c | c | c | c | c |}  
 \hline
 $m$ & $R(\phi)$ & $R(\phi^\prime)$ &                     $R(\phi^{\prime\prime})$ & $R(\phi^{(3)})$  &                     $R(\phi^{(4)})$ & Time (s)\\  
 \hline
   70  &      0.63E-11  &      0.66E-09  &      0.67E-09  &      0.45E-09  &      0.45E-09 
          &      0.32E+01 \\
   80  &      0.10E-16  &      0.10E-14  &      0.10E-14  &      0.82E-15  &      0.82E-15
          &      0.33E+01  \\
   90  &      0.11E-23  &      0.11E-21  &      0.11E-21  &      0.99E-22  &      0.99E-22
          &      0.37E+01  \\
  100  &      0.12E-25  &      0.30E-27  &      0.23E-29  &      0.30E-29  &      0.30E-29
          &      0.40E+01  \\
  110  &      0.22E-25  &      0.22E-26  &      0.56E-30  &      0.61E-30  &      0.61E-30
          &      0.43E+01  \\
  120  &      0.67E-25  &      0.34E-26  &      0.11E-30  &      0.63E-30  &      0.63E-30
          &      0.48E+01  \\
 \hline
 \end{tabular}
 \end{table} 
 \begin{table} 
 \caption{Relative residual sizes for                    boundary layer}
 \label{residuals:bdlayer}.
 \small 
 \center 
 \begin{tabular}{|l |c |c |c |c |c |c |c | } 
 \hline
& \multicolumn{    7             }{c|}{Iteration number}\\   \cline{2-    8}
 $m$ &   1 &   2 &   3 &   4 &   5 &   6 &   7  \\
 \hline
   70  &      0.51E+02  &      0.97E-05  &      0.33E-07  &      0.63E-14  &      0.22E-16  &      0.17E-22  &      0.86E-23    \\  
   80 &     0.90E-04 &     0.15E-10 &     0.91E-19 &     0.11E-22 &  &    &     \\ 
   90 &     0.10E-10 &     0.17E-17 &     0.21E-22 &     0.93E-23 &  &   &    \\ 
  100 &     0.24E-18 &     0.37E-21 &     0.83E-23 &     &    &     &    \\ 
  110 &     0.14E-19 &     0.96E-22 &     0.80E-23 & &    &    &      \\ 
  120 &     0.22E-19 &     0.25E-21 &     0.12E-22 &      &      &    &    \\ 
 \hline
 \end{tabular}
 \end{table}

\end{document}